\documentclass{article}

\usepackage[english]{babel}
\usepackage[utf8]{inputenc}
\usepackage[letterpaper, hmargin=1.4in, vmargin=1in]{geometry}
\usepackage{amsmath, amsthm, amssymb, dsfont}
\usepackage{graphicx}
\usepackage[margin=40pt,font=small,labelfont=bf]{caption}
\usepackage[textsize=scriptsize,colorinlistoftodos]{todonotes}
\usepackage{hyperref}
\hypersetup{
    colorlinks,
    linkcolor={red!60!black},
    citecolor={blue!60!black},
    urlcolor={blue!80!black}
}
\usepackage[capitalize, noabbrev]{cleveref}
\usepackage{comment}
\usepackage{booktabs}
\usepackage{tikz}
\usetikzlibrary{decorations.pathreplacing}
\usepackage{pgfplots}
\pgfplotsset{compat=1.16}

\theoremstyle{definition}

\newtheorem*{definition}{Definition}

\begin{document}

\title{\textbf{Magic Triangles}}

\author{Gabriel Hale, Bjorn Vogen, Matthew Wright\footnote{St.\ Olaf College, Northfield, MN}}

\date{June 2021}

\maketitle

\begin{abstract}
Magic squares are well-known arrangements of integers with common row, column, and diagonal sums. 
Various other magic shapes have been proposed, but triangles have been somewhat overlooked.
We introduce certain triangular arrangements of integers with common sums in three directions, which we call magic triangles.
For small sizes of these triangles, we count the number of unique magic triangles and examine distributions of integers at different positions within them. 
While we cannot enumerate the number of magic triangles at larger sizes, we offer a simulated annealing method for finding magic triangles.
\end{abstract}

\section{Introduction}

Consider an equilateral triangle, partitioned into sixteen congruent subtriangles as in the following diagram. 
Each of the numbers $1, 2, 3, \ldots, 16$ is placed into one of the subtriangles as shown.

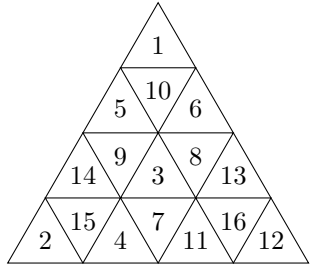
\begin{center}
  \begin{tikzpicture}
    \draw (0,0) -- (4,0) -- (2,3.464) -- cycle;
    \draw (1,1.732) -- (2,0) -- (3,1.732) -- cycle;
    \draw (0.5,0.866) -- (3.5,0.866) -- (3,0) -- (1.5,2.598) -- (2.5,2.598) -- (1,0) -- cycle;
    
    \node at (2,2.9) {$1$};
    
    \node at (1.5,2.05) {$5$};
    \node at (2,2.3) {$10$};
    \node at (2.5,2.05) {$6$};
    
    \node at (1,1.17) {$14$};
    \node at (1.5,1.42) {$9$};
    \node at (2,1.17) {$3$};
    \node at (2.5,1.42) {$8$};
    \node at (3,1.17) {$13$};
    
    \node at (0.5,0.3) {$2$};
    \node at (1,0.55) {$15$};
    \node at (1.5,0.3) {$4$};
    \node at (2,0.55) {$7$};
    \node at (2.5,0.3) {$11$};
    \node at (3,0.55) {$16$};
    \node at (3.5,0.3) {$12$};
  \end{tikzpicture}
\end{center}

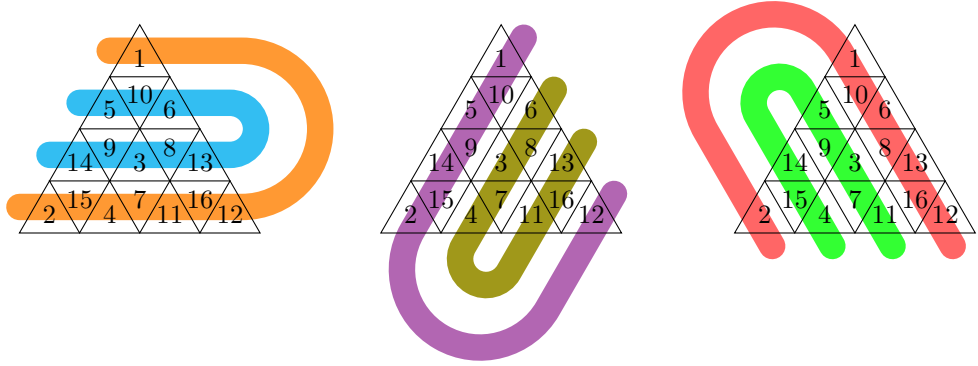
\begin{figure}[h!]
  \centering
  \begin{tikzpicture}[scale=0.8,baseline=0pt]
    \draw[cyan, line width=10pt, line cap=round, opacity=0.8] (0.5,1.299) -- ++(3,0) arc (-90:90:0.433) -- (1,2.165);
    \draw[orange, line width=10pt, line cap=round, opacity=0.8] (0,0.433) -- ++(3.7,0) arc (-90:90:1.299) -- (1.5,3.031);
    
    \draw (0,0) -- (4,0) -- (2,3.464) -- cycle;
    \draw (1,1.732) -- (2,0) -- (3,1.732) -- cycle;
    \draw (0.5,0.866) -- (3.5,0.866) -- (3,0) -- (1.5,2.598) -- (2.5,2.598) -- (1,0) -- cycle;
    
    \node at (2,2.9) {$1$};
    
    \node at (1.5,2.05) {$5$};
    \node at (2,2.3) {$10$};
    \node at (2.5,2.05) {$6$};
    
    \node at (1,1.17) {$14$};
    \node at (1.5,1.42) {$9$};
    \node at (2,1.17) {$3$};
    \node at (2.5,1.42) {$8$};
    \node at (3,1.17) {$13$};
    
    \node at (0.5,0.3) {$2$};
    \node at (1,0.55) {$15$};
    \node at (1.5,0.3) {$4$};
    \node at (2,0.55) {$7$};
    \node at (2.5,0.3) {$11$};
    \node at (3,0.55) {$16$};
    \node at (3.5,0.3) {$12$};
  \end{tikzpicture}
  \hspace{10pt}
  \begin{tikzpicture}[scale=0.8,baseline=0pt]
    \draw[olive, line width=10pt, line cap=round, opacity=0.9] (2.875,2.382) -- ++(240:3) arc (150:330:0.433) -- (3.375,1.516);
    \draw[violet, line width=10pt, line cap=round, opacity=0.6] (2.375,3.248) -- ++(240:3.7) arc (150:330:1.299) -- (3.875,0.650);
    
    \draw (0,0) -- (4,0) -- (2,3.464) -- cycle;
    \draw (1,1.732) -- (2,0) -- (3,1.732) -- cycle;
    \draw (0.5,0.866) -- (3.5,0.866) -- (3,0) -- (1.5,2.598) -- (2.5,2.598) -- (1,0) -- cycle;
    
    \node at (2,2.9) {$1$};
    
    \node at (1.5,2.05) {$5$};
    \node at (2,2.3) {$10$};
    \node at (2.5,2.05) {$6$};
    
    \node at (1,1.17) {$14$};
    \node at (1.5,1.42) {$9$};
    \node at (2,1.17) {$3$};
    \node at (2.5,1.42) {$8$};
    \node at (3,1.17) {$13$};
    
    \node at (0.5,0.3) {$2$};
    \node at (1,0.55) {$15$};
    \node at (1.5,0.3) {$4$};
    \node at (2,0.55) {$7$};
    \node at (2.5,0.3) {$11$};
    \node at (3,0.55) {$16$};
    \node at (3.5,0.3) {$12$};
  \end{tikzpicture}
  \hspace{10pt}
  \begin{tikzpicture}[scale=0.8,baseline=0pt]
    \draw[green, line width=10pt, line cap=round, opacity=0.8] (2.625,-0.217) -- ++(120:3) arc (30:210:0.433) -- (1.625,-0.217);
    \draw[red, line width=10pt, line cap=round, opacity=0.6] (3.625,-0.217) -- ++(120:3.7) arc (30:210:1.299) -- (0.625,-0.217);
  
    \draw (0,0) -- (4,0) -- (2,3.464) -- cycle;
    \draw (1,1.732) -- (2,0) -- (3,1.732) -- cycle;
    \draw (0.5,0.866) -- (3.5,0.866) -- (3,0) -- (1.5,2.598) -- (2.5,2.598) -- (1,0) -- cycle;
    
    \node at (2,2.9) {$1$};
    
    \node at (1.5,2.05) {$5$};
    \node at (2,2.3) {$10$};
    \node at (2.5,2.05) {$6$};
    
    \node at (1,1.17) {$14$};
    \node at (1.5,1.42) {$9$};
    \node at (2,1.17) {$3$};
    \node at (2.5,1.42) {$8$};
    \node at (3,1.17) {$13$};
    
    \node at (0.5,0.3) {$2$};
    \node at (1,0.55) {$15$};
    \node at (1.5,0.3) {$4$};
    \node at (2,0.55) {$7$};
    \node at (2.5,0.3) {$11$};
    \node at (3,0.55) {$16$};
    \node at (3.5,0.3) {$12$};
  \end{tikzpicture}
  \caption{Horizontal, positive slope, and negative slope sums.}
  \label{sumfig}
\end{figure}

Examining the sums along rows and diagonals reveals an interesting property. 
First, consider the sums along horizontal rows, highlighted in \cref{sumfig} (left).
The sum of the entries in the top and bottom rows is $1 + (2 + 15 + 4 + 7 + 11 + 16 + 12) = 68$. 
The sum of the entries in the middle two rows is $(5 + 10 + 6) + (14 + 9 + 3 + 8 + 13) = 68$. 
We see that we have obtained the same sum twice.

Next, consider the sums along diagonals with positive slope, illustrated in \cref{sumfig} (middle).
The sum of the entries in the right corner and the leftmost diagonal is $12 + (1 + 10 + 5 + 9 + 14 + 15 + 2) = 68$.
Adding up the entries in the middle two diagonals, we find the same sum again: $(4 + 7 + 3 +8 + 6) + (11 + 16 + 13) = 68$.

Finally, compute the sums along the diagonals of negative slope, illustrated in \cref{sumfig} (right).
The sum of the entries in the right corner and the leftmost diagonal is $2 + (1 + 10 + 6 + 8 + 13 + 16 + 12) = 68$, and the sum of the middle two diagonals is $(14 + 15 + 4) + (5 + 9 + 3 + 7 + 11) = 68$.

Fascinatingly, we have obtained the same sum, $68$, by adding up the entries in six pairs of rows or diagonals in this triangular arrangement of integers. 
We call this arrangement a \emph{magic triangle}. 
Several questions are immediate: Is this the only magic triangle? If not, how many other magic triangles are there? Can triangular arrangements of other sizes be magic?

While magic squares and other magic shapes have been studied extensively, we are not aware of any literature dealing with the triangular arrangement and sums described above.
Among the three regular polygons that tile the plane, magic squares are well known and there is a unique magic hexagon \cite{Trigg}, but magic triangles have scarcely been considered.
Indeed, the term ``magic triangle'' seems to refer most commonly to \emph{perimeter magic triangles}---arrangements of integers along the perimeter of a triangle, such that the sums along each side of the triangle are the same \cite{pmt}.
For instance, in his classic ``Magic Designs'' article, Robert Ely discusses certain triangular arrangements, but does not suggest the sums that we study in this article  \cite{Ely}.
Nor does our concept appear in Clifford Pickover's encyclopedic work on magic squares, circles, and stars \cite{Pickover}.
Our magic triangles fill a gap in the literature on magic shapes.

\section{Magic triangle definition}

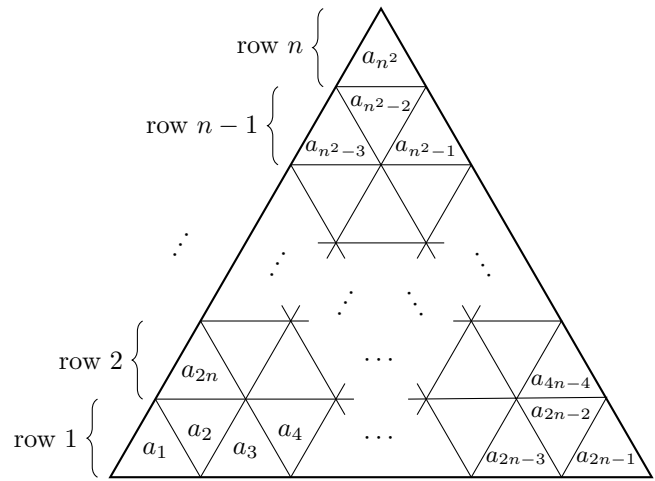
\begin{figure}[h!]
  \centering
  \begin{tikzpicture}[scale=1.2]
    \draw[thick] (0,0) -- (6,0) -- (3,5.196) -- cycle;
    
    \draw (0.5,0.866) -- ++(2.2,0);
    \draw (1,1.732) -- ++(1.2,0);
    \draw (0.5,0.866) -- (1,0) -- ++(60:2.2);
    \draw (1,1.732) -- (2,0) -- ++(60:1.2);
    \draw (2,1.732) -- ++(120:0.2);
    \draw (2,1.732) -- ++(300:1.2);
    
    \draw (3.3,0.866) -- (5.5,0.886) -- (5,0) -- ++(120:2.2);
    \draw (3.8,1.732) -- (5,1.732) -- (4,0) -- ++(120:1.2);
    \draw (4,1.732) -- ++(60:0.2);
    \draw (4,1.732) -- ++(240:1.2);
    
    \draw (2.5,4.330) -- (3.5,4.330) -- ++(-120:2.2);
    \draw (2.5,4.330) -- ++(-60:2.2);
    \draw (2,3.464) -- (4,3.464);
    \draw (2.3,2.598) -- ++(1.4,0);
    \draw (2,3.464) -- ++(-60:1.2);
    \draw (4,3.464) -- ++(-120:1.2);
    
    \node at (3,0.43) {$\cdots$};
    \node at (3,1.3) {$\cdots$};
    \node[rotate=60] at (1.875,2.382) {$\cdots$};
    \node[rotate=60] at (2.625,1.949) {$\cdots$};
    \node[rotate=120] at (4.125,2.382) {$\cdots$};
    \node[rotate=120] at (3.375,1.949) {$\cdots$};
    
    %\draw (0.5,0.866) -- (4.5,0.866) -- (4,0) -- (2,3.464) -- (3,3.464) -- (1,0) -- cycle;
    %\draw (1,1.732) -- (4,1.732) -- (3,0) -- (1.5,2.598) -- (3.5,2.598) -- (2,0) -- cycle;
    
    \node at (0.5,0.3) {$a_1$};
    \node at (1,0.5) {$a_2$};
    \node at (1.5,0.3) {$a_3$};
    \node at (2,0.5) {$a_4$};
    \node at (4.5,0.2) {\small{$a_{2n-3}$}};
    \node at (5,0.7) {\small{$a_{2n-2}$}};
    \node at (5.5,0.2) {\small{$a_{2n-1}$}};
    \node at (1,1.15) {$a_{2n}$};
    \node at (5,1.05) {\small{$a_{4n-4}$}};
    
    \node at (2.5,3.63) {\small{$a_{n^2-3}$}};
    \node at (3,4.13) {\small{$a_{n^2-2}$}};
    \node at (3.5,3.63) {\small{$a_{n^2-1}$}};
    
    \node at (3,4.63) {$a_{n^2}$};
    
    \draw [decorate,decoration={brace,amplitude=4pt},xshift=-4pt,yshift=0pt] (0,0) -- ++(0,0.866) node [midway, left,xshift=-4pt] {row 1};
    \draw [decorate,decoration={brace,amplitude=4pt},xshift=-4pt,yshift=0pt] (0.5,0.866) -- ++(0,0.866) node [midway, left,xshift=-4pt] {row 2};
    \node[rotate=60] at (0.8,2.6) {$\cdots$};
    \draw [decorate,decoration={brace,amplitude=4pt},xshift=-4pt,yshift=0pt] (2,3.464) -- ++(0,0.866) node [midway, left,xshift=-4pt] {row $n-1$};
    \draw [decorate,decoration={brace,amplitude=4pt},xshift=-4pt,yshift=0pt] (2.5,4.33) -- ++(0,0.866) node [midway, left,xshift=-4pt] {row $n$};
  \end{tikzpicture}
  \caption{Entries and rows of an $n$-level triangle.}
  \label{nlevel}
\end{figure}

To define a \emph{magic triangle}, we formalize the example above.
Let $n$ be a positive integer and consider an equilateral triangle partitioned into $n^2$ subtriangles, as shown in \cref{nlevel}.
Place the integers $1, 2, \ldots, n^2$ into the subtriangles.
We label these entries $a_1, a_2, \ldots, a_{n^2}$ as in the figure, choosing (somewhat arbitrarily) to index the entries from left-to-right, bottom-to-top.
We emphasize that $(a_1, a_2, \ldots, a_{n^2})$ may be any permutation of $1, 2, \ldots, n^2$.

Observe that there are $n$ rows of subtriangles, as labeled in \cref{nlevel}. 
For an integer $k \in \{1, \ldots, n\}$, let $r_k$ be the sum of the entries in row $k$. 
%Specifically, \[ r_k = \sum_{i = 2nk-k^2+2k-2n}^{2nk-k^2} a_i. \]
We then define the \emph{horizontal sums} $h_k = r_k + r_{n-k}$.
By symmetry, $h_k = h_{n-k}$.
If $n$ is odd, then the ``middle'' horizontal sum is twice the sum of the entries in the middle row of the large triangle: $h_{\frac{n+1}{2}} = 2 r_{\frac{n+1}{2}}$.
All other horizontal sums are sums of $n$ distinct integers in the large triangle.

Adding up all of the horizontal sums, $h_1 + \cdots + h_n$, we add each entry in the triangle exactly twice. 
That is,
\[ \sum_{k=1}^n h_k = 2 \sum_{i=1}^{n^2} a_i = 2 \sum_{i=1}^{n^2} i = 2\cdot \frac{n^2(n^2+1)}{2} = n^2(n^2+1). \]
Thus, if all of the $n$ horizontal sums are to be equal, it must be that they are each equal to $n(n^2+1)$.
Since $h_k = h_{n-k}$, it suffices to require that
\[ h_1 = h_2 = \cdots = h_{\lceil \frac{n}{2} \rceil} = n(n^2+1), \]
where $\lceil \frac{n}{2} \rceil$ denotes the smallest integer greater than or equal to $\frac{n}{2}$.

We similarly define \emph{positive slope sums} $p_1, \ldots, p_n$ and \emph{negative slope sums} $q_1, \ldots, q_n$.
Rather than defining these sums in terms of the indexes in \cref{nlevel}, we approach them by rotating the triangular arrangement $120^\circ$.
Specifically, rotating the triangular arrangement $120^\circ$ counterclockwise makes the positive slope diagonals become horizontal; the $p_k$ can be defined as the horizontal sums of this rotated arrangement.
Likewise, rotating the original triangular arrangement $120^\circ$ clockwise makes the negative slope diagonals become horizontal; the $q_k$ can be defined as the horizontal sums of this rotated arrangement.
The positive and negative slope sums generalize the sums illustrated in the middle and right diagrams, respectively, in \cref{sumfig}.

This leads to our formal definition of a magic triangle:
\begin{definition}
  A triangular arrangement of $1, 2, \ldots, n^2$ is an \emph{$n$-level magic triangle} if the horizontal, positive slope, and negative slope sums are all equal. 
  That is, if $h_k = p_k = q_k = n(n^2+1)$ for all $k \in \left\{1, 2, \ldots, \lceil \frac{n}{2} \rceil \right\}$.
\end{definition}

\section{How many magic triangles?}

How many $n$-level magic triangles are there, for each positive integer $n$? While this is a difficult question for large $n$, we provide answers for a few small $n$.

We count magic triangles \emph{up to symmetry}, meaning that we do not distinguish between rotations and reflections of a given triangular arrangement.
Specifically, any triangular arrangement can be flipped across each of its altitudes and rotated $120^\circ$ clockwise or counterclockwise, producing six arrangements that we regard as the same up to symmetry.
Let $T_n$ denote the number of $n$-level magic triangles, up to symmetry.

First, the $n=1$ case is trivial: there is only one triangle containing the integer $1$, which is magic, so $T_1=1$.

\begin{figure}[h!]
    \centering
    \begin{tikzpicture}[scale=0.9]
      \draw[thick] (0,0) -- (2,0) -- (1,1.732) -- cycle;
      \draw (0.5,0.866) -- (1,0) -- (1.5,0.866) -- cycle;
      \node at (0.5,0.3) {$a_1$};
      \node at (1,0.5) {$a_2$};
      \node at (1.5,0.3) {$a_3$};
      \node at (1,1.15) {$a_4$};
      
      \draw (3,0) -- ++(0,1.7);
    \end{tikzpicture}
    \hspace{18pt}
    \begin{tikzpicture}[scale=0.9]
      \draw[thick] (0,0) -- (2,0) -- (1,1.732) -- cycle;
      \draw (0.5,0.866) -- (1,0) -- (1.5,0.866) -- cycle;
      \node at (0.5,0.3) {$1$};
      \node at (1,0.5) {$2$};
      \node at (1.5,0.3) {$3$};
      \node at (1,1.15) {$4$};
    \end{tikzpicture}
    \hspace{15pt}
    \begin{tikzpicture}[scale=0.9]
      \draw[thick] (0,0) -- (2,0) -- (1,1.732) -- cycle;
      \draw (0.5,0.866) -- (1,0) -- (1.5,0.866) -- cycle;
      \node at (0.5,0.3) {$1$};
      \node at (1,0.5) {$3$};
      \node at (1.5,0.3) {$2$};
      \node at (1,1.15) {$4$};
    \end{tikzpicture}
    \hspace{15pt}
    \begin{tikzpicture}[scale=0.9]
      \draw[thick] (0,0) -- (2,0) -- (1,1.732) -- cycle;
      \draw (0.5,0.866) -- (1,0) -- (1.5,0.866) -- cycle;
      \node at (0.5,0.3) {$1$};
      \node at (1,0.5) {$4$};
      \node at (1.5,0.3) {$2$};
      \node at (1,1.15) {$3$};
    \end{tikzpicture}
    \hspace{15pt}
    \begin{tikzpicture}[scale=0.9]
      \draw[thick] (0,0) -- (2,0) -- (1,1.732) -- cycle;
      \draw (0.5,0.866) -- (1,0) -- (1.5,0.866) -- cycle;
      \node at (0.5,0.3) {$2$};
      \node at (1,0.5) {$1$};
      \node at (1.5,0.3) {$3$};
      \node at (1,1.15) {$4$};
    \end{tikzpicture}
    \caption{The four entries of a $2$-level triangle, and the four $2$-level magic triangles, up to symmetry.}
    \label{twolevel}
\end{figure}
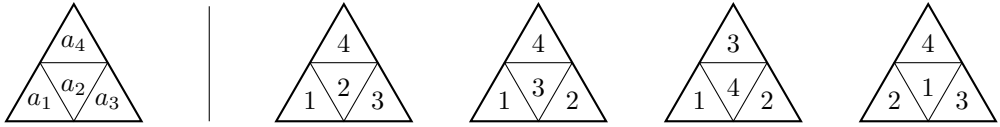

The $n=2$ case provides a simple illustration of symmetry.
A $2$-level triangle has only one horizontal sum, one positive slope sum, and one negative slope sum.
Thus, each of the integers $1, 2, 3, 4$ are a part of each sum, so clearly any arrangement of these integers yields a magic triangle.
Labeling the four entries as in \cref{twolevel} (left), we can count up to symmetry by fixing the relative order of the corner entries, such as $a_1 < a_3 < a_4$.
With this constraint, we find the four arrangements shown in \cref{twolevel} (right), and thus $T_2 = 4$.

Alternatively, there are $4!=24$ unconstrained permutations of the four entries in the $2$-level triangle. Dividing by the number of symmetries of a triangle gives $T_2 = \frac{4!}{6} = 4$.

\subsection*{3-level triangles}

The $n=3$ case is the smallest case in which not every triangular arrangement is magic.
While we could consider all $9!$ arrangements in an exhaustive search, the geometry of the triangle allows us to reduce the work required to count $3$-level magic triangles.  

A $3$-level magic triangle must satisfy two sums in each direction (horizontal, positive slope, and negative slope). For instance, the horizontal sums must satisfy $h_1 = h_2 = 30$. 
However, it suffices to focus on the middle row each direction:
if we guarantee that the middle row sum is $a_6 + a_7 + a_8 = 15$, then $h_2 = 30$ because $h_2$ is twice the sum of the middle row. 
Since the sum of all entries in the triangle is $45$, if $h_2=30$ then it follows that $h_1 = 30$ as well.

\begin{figure}[h!]
  \centering
  \begin{tikzpicture}[scale=1.2,baseline=0pt]
    \draw[thick] (0,0) -- (3,0) -- (1.5,2.598) -- cycle;
    \draw (0.5,0.866) -- (1,0) -- (2,1.732) -- (1,1.732) -- (2,0) -- (2.5,0.866) -- cycle;
    
    \node at (0.5,0.3) {$a_1$};
    \node at (1,0.5) {$a_2$};
    \node at (1.5,0.3) {$a_3$};
    \node at (2,0.5) {$a_4$};
    \node at (2.5,0.3) {$a_5$};
    \node at (1,1.15) {$a_6$};
    \node at (1.5,1.35) {$a_7$};
    \node at (2,1.15) {$a_8$};
    \node at (1.5,2) {$a_9$};
  \end{tikzpicture}
  \hspace{30pt}
  \begin{tikzpicture}[scale=1.2,baseline=0pt]
    \draw[cyan, line width=12pt, line cap=round, opacity=0.8] (0.5,1.299) -- ++(2,0);
    \draw[violet, line width=12pt, line cap=round, opacity=0.7] (2.375,1.516) -- ++(240:2);
    \draw[green, line width=12pt, line cap=round, opacity=0.8] (1.625,-0.217) -- ++(120:2);
    
    \draw (1,0) -- (2,0) -- (2.5,0.866) -- (2,1.732) -- (1,1.732) -- (0.5,0.866) -- cycle;
    \draw (0.5,0.866) -- (2.5,0.866);
    \draw (1,0) -- (2,1.732);
    \draw (2,0) -- (1,1.732);
    
    \draw[darkgray] (0.4,0.692) -- (0,0) -- (0.8,0);
    \draw[darkgray] (2.6,0.692) -- (3,0) -- (2.2,0);
    \draw[darkgray] (1.9,1.902) -- (1.5,2.598) -- (1.1,1.902);
    
  \end{tikzpicture}
  \caption{A 3-level triangle and its central hexagon with three sums highlighted.}
  \label{threelevel}
\end{figure}
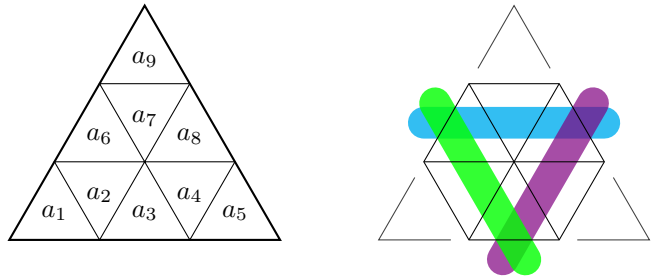

Thus, we want the sum of the entries in the middle row of the triangle, in each direction, to equal $15$.
This leads us to focus on the ``central hexagon''---the six non-corner entries of the triangle, displayed at right in \cref{threelevel}.
The three middle rows (highlighted in \cref{threelevel}) are comprised of entries in this central hexagon.
We have a magic triangle exactly when
\begin{equation}\label{level3sums}
  a_6 + a_7 + a_8 = a_3 + a_4 + a_8 = a_2 + a_3 + a_6 = 15. 
\end{equation}

The three sums in \cref{level3sums} impose the symmetries of a triangle on to the central hexagon, illustrated by the highlighed triangle in \cref{threelevel} (right).
We count solutions up to these symmetries by imposing the additional requirement $a_3 < a_6 < a_8$.
By exhaustively checking sets of six unique integers from $\{1, \ldots, 9\}$, we find $16$ solutions.
Each of the $16$ solutions can be flipped and rotated $6$ ways within the central hexagon, for a total of $96$ solutions to \cref{level3sums}.

Any solution to \cref{level3sums} extends to a magic triangle by placing the three integers not in the central hexagon into the corner entries of the triangle.
By requiring the ordering $a_1 < a_5 < a_9$ on the corners of the triangle, we have a bijection between $3$-level magic triangles up to symmetry and solutions to \cref{level3sums}.
Therefore, $T_3 = 96$.

\subsection*{4-level triangles}

The $n=4$ case involves $16$ integers that can be arranged in $16! \approx 2 \times 10^{13}$ ways, which is too many for us to consider exhaustively.
As in the previous case, the geometry of the triangle allows us to reduce the search space and determine $T_4$.

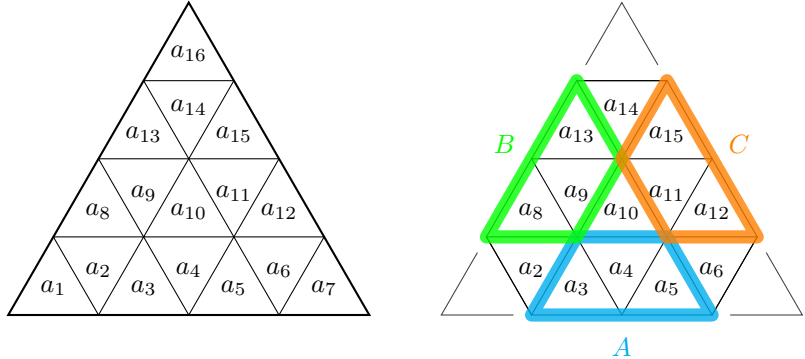
\begin{figure}[h!]
  \centering
  \begin{tikzpicture}[scale=1.2,baseline=0pt]
    \draw[thick] (0,0) -- (4,0) -- (2,3.464) -- cycle;
    \draw (0.5,0.866) -- (1,0) -- (2.5,2.598) -- (1.5,2.598) -- (3,0) -- (3.5,0.866) -- cycle;
    \draw (1,1.732) -- (3,1.732) -- (2,0) -- cycle;
    
    \node at (0.5,0.3) {$a_1$};
    \node at (1,0.5) {$a_2$};
    \node at (1.5,0.3) {$a_3$};
    \node at (2,0.5) {$a_4$};
    \node at (2.5,0.3) {$a_5$};
    \node at (3,0.5) {$a_6$};
    \node at (3.5,0.3) {$a_7$};
    \node at (1,1.15) {$a_8$};
    \node at (1.5,1.35) {$a_9$};
    \node at (2,1.15) {$a_{10}$};
    \node at (2.5,1.35) {$a_{11}$};
    \node at (3,1.15) {$a_{12}$};
    \node at (1.5,2) {$a_{13}$};
    \node at (2,2.3) {$a_{14}$};
    \node at (2.5,2) {$a_{15}$};
    \node at (2,2.9) {$a_{16}$};
  \end{tikzpicture}
  \hspace{20pt}
    \begin{tikzpicture}[scale=1.2,baseline=0pt]
    \draw (1,0) -- (3,0) -- (3.5,0.866) -- (2.5,2.598) -- (1.5,2.598) -- (0.5,0.866) -- cycle;
    \draw (0.5,0.866) -- (1,0) -- (2.5,2.598) -- (1.5,2.598) -- (3,0) -- (3.5,0.866) -- cycle;
    \draw (1,1.732) -- (3,1.732) -- (2,0) -- cycle;
    
    \draw[darkgray] (0.4,0.692) -- (0,0) -- (0.8,0);
    \draw[darkgray] (3.6,0.692) -- (4,0) -- (3.2,0);
    \draw[darkgray] (1.6,2.772) -- (2,3.464) -- (2.4,2.772);
    
    \draw[cyan, line width=5pt, line join=round, opacity=0.8] (1,0) -- (3,0) -- (2.5,0.866) -- (1.5,0.866) -- cycle;
    \draw[green, line width=5pt, line join=round, opacity=0.8] (0.5,0.866) -- (1.5,0.866) -- (2,1.732) -- (1.5,2.598) -- cycle;
    \draw[orange, line width=5pt, line join=round, opacity=0.8] (3.5,0.866) -- (2.5,0.866) -- (2,1.732) -- (2.5,2.598) -- cycle;
    
    \node[cyan] at (2,-0.35) {$A$};
    \node[green] at (0.7,1.9) {$B$};
    \node[orange] at (3.3,1.9) {$C$};
    
    \node at (1,0.5) {$a_2$};
    \node at (1.5,0.3) {$a_3$};
    \node at (2,0.5) {$a_4$};
    \node at (2.5,0.3) {$a_5$};
    \node at (3,0.5) {$a_6$};
    \node at (1,1.15) {$a_8$};
    \node at (1.5,1.35) {$a_9$};
    \node at (2,1.15) {$a_{10}$};
    \node at (2.5,1.35) {$a_{11}$};
    \node at (3,1.15) {$a_{12}$};
    \node at (1.5,2) {$a_{13}$};
    \node at (2,2.3) {$a_{14}$};
    \node at (2.5,2) {$a_{15}$};
  \end{tikzpicture}
  \caption{A 4-level triangle (left) and its central hexagon (right) with sets $A$, $B$, and $C$ highlighted.}
  \label{fourlevel}
\end{figure}

As before, we temporarily ignore the corner entries and focus on the central hexagon.
The central hexagon is now not a regular hexagon, but consists of the $13$ entries shown at right in \cref{fourlevel}. 
We seek to assign integers to the central hexagon entries in a way that satisfies the magic triangle sums $h_2 = p_2 = q_2 = 68$, but we impose additional requirements to limit the number of possibilities that we must check.

Within the central hexagon, we notice that three groups of three entries always appear together in the sums that define the magic triangle.
For example, $a_3, a_4, a_5$ occur together in the outer horizontal sum $h_1$, and also in the inner sums of positive and negative slope, $p_2$ and $q_2$. Thus, permutations among these three entries do not change the sums.
We let $A = \{a_3, a_4, a_5\}$, as illustrated in \cref{fourlevel}.
Similarly, $a_8, a_9, a_{13}$ always appear together in the sums, so we let $B = \{a_8, a_9, a_{13}\}$.
We let $C = \{a_{11}, a_{12}, a_{15}\}$, the third set of three entries that always appear together in the sums. 
As we search for solutions, we may thus impose a relative order on the entries within each group:
\begin{equation}\label{level4grouporders}
    a_3 < a_4 < a_5, \qquad a_8 < a_9 < a_{13}, \qquad\text{and}\qquad a_{11} < a_{12} < a_{15}. 
\end{equation}

Furthermore, the central hexagon itself has six symmetries.
Rotating or reflecting the central hexagon within the $4$-level triangle simply permutes the sums that define the magic triangle.
We thus impose the additional requirement
\begin{equation}\label{level4hexorder}
    a_2 < a_6 < a_{14}.
\end{equation}

Thus, we seek an assignment of $13$ integers from $\{1, \ldots, 16\}$ to the entries of the central hexagon such that $h_2 = p_2 = q_3 = 68$ and Inequalities \eqref{level4grouporders} and \eqref{level4hexorder} hold.
The inequalities reduce the search space to 
\[ \binom{16}{13}\cdot \frac{13!}{6^4}= 2{,}690{,}688{,}000 \]
possible assignments of distinct integers to the central hexagon entries; for each we check whether $h_2 = p_2 = q_3 = 68$ holds.
We find $184{,}056$ central hexagon solutions.

Each central hexagon solution yields $6^4$ magic triangles: the central hexagon may be positioned in the triangle in $6$ ways, and the entries within each region $A$, $B$, and $C$ may be permuted $6$ ways. 
The corner entries of the triangle are then determined: the three integers not in the central hexagon are the corner entries of the magic triangle, with $a_1 < a_7 < a_{16}$ to account for the symmetries of the triangle.
Therefore, the number of $4$-level magic triangles, up to symmetry, is
\[ T_4 = 184{,}056 \cdot 6^4 = 238{,}536{,}576. \]

\subsection*{More than $4$ levels}

Our strategy for counting $3$- and $4$-level magic triangles generalizes to triangles with $5$ or more levels, but the number of possible assignments of integers becomes too large for us to check exhaustively, even when using geometry to limit the search. 
Thus, we are unable to give an exact value for $T_n$ with $n \ge 5$.
This mirrors the situation for magic squares, in that the number of $m \times m$ magic squares is not known for $m \ge 6$ \cite{magicSquareSeq}.

We generated ten billion random arrangements of entries in a $5$-level triangle and found not even a single magic triangle. This suggests that the frequency of  magic triangles among all $5$-level triangles is less than one in $10^{10}$. 
Thus, it seems that there are are fewer than 
\[ \frac{25!}{6}\cdot \frac{1}{10^{10}} \approx 2.6 \times 10^{14} \]
5-level magic triangles (up to symmetry), and possibly many fewer.

Instead of random generation, an effective method of finding magic triangles is \emph{simulated annealing}, a probabilistic algorithm that shuffles entries in a way that is likely to move from a random triangular arrangement to a magic triangle.
Simulated annealing is often used to solve combinatorial optimization problems \cite{ShoMen}, and is similar to the method used by Kitajima and Kikuchi to count rare cases of magic squares \cite{Numerous}.

The simulated annealing algorithm starts with a random triangular arrangement of integers and repeatedly performs the following step: two entries are selected at random and possibly swapped.
If swapping the entries moves the arrangement closer to a magic triangle, in the sense that the magic triangle sums move closer to the target value, then the swap is performed.
Otherwise, the swap is performed with a probability that decreases with each step, allowing the algorithm to avoid getting trapped in a configuration from which no swap will result in a magic triangle.
The algorithms stops when a magic triangle is found.

\begin{figure}[h!]
    \centering

    \begin{tikzpicture}
        \begin{axis}[width=3.5in, height=2.4in, xlabel={$n$}, xtick={3,4,...,8}, ylabel={number of steps}, ymode=log, ytick={10,1e2,1e3,1e4,1e5,1e6}, ymin=20, yminorticks=true, grid=major, legend style={at={(0.95,0.35)}}]

            % mean values
            \addplot[color=red,mark=x] coordinates { (3,94) (4, 277) (5, 5933) (6, 7696) (7, 45637) (8, 246845) };
            
            % median values
            \addplot[color=blue,mark=*] coordinates { (3,71) (4, 205) (5, 3179) (6, 5536) (7, 33309) (8, 174527) };
            
            \legend{mean, median}
        \end{axis}
    \end{tikzpicture}
    
    \caption{Mean and median number of steps required to find magic triangles for $n$ from $3$ to $8$, obtained from $10{,}000$ simulated annealing experiments for each $n$.}
    \label{fig:sim_annel}
\end{figure}
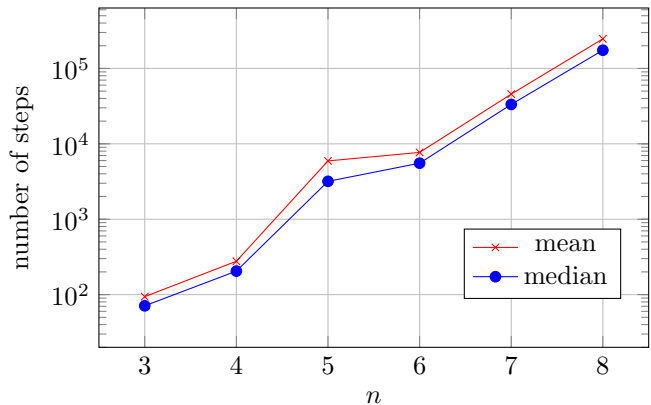

Using simulated annealing, we have found magic triangles up to size $n=10$.
In general, larger $n$ requires more steps before a magic triangle is found.
The simulated annealing algorithm also involves several parameters that must be experimentally tuned for each size of triangle.
After tuning the parameters for each $n$ from $3$ to $8$, we ran $10{,}000$ simulations for each $n$. For each simulation, we recorded the number of steps required to find a magic triangle.
\Cref{fig:sim_annel} displays the mean and median number of steps for each $n$.

For each $n$, the distribution of the number of steps is skewed right: the algorithm often finds magic triangles quickly, but sometimes requires a large number of steps. As such, the median is perhaps the more useful measure of the algorithm's efficiency: in half of the simulations, the algorithm finds a magic triangle in at most the median number of steps.

\section{Distributions of integers in different positions}

Within all magic triangles of a fixed size, are certain integers more commonly found in certain positions? 
We discover interesting patterns in the frequencies with which integers appear in different positions within $3$-level and $4$-level magic triangles.

\subsection*{3-level triangles}

The $3$-level triangle contains three orbits of entries under rotations and reflections of the triangle. 
Referring to \cref{threelevel}, we call $\{a_1, a_5, a_9\}$ the \emph{corner entries}, $\{a_3, a_6, a_8\}$ the \emph{border entries}, and $\{a_2, a_4, a_7\}$ the \emph{interior entries}.
We counted how many times each integer $1, 2, \ldots, 9$ appears in the corner, border, and interior entries among all $96$ magic triangles.

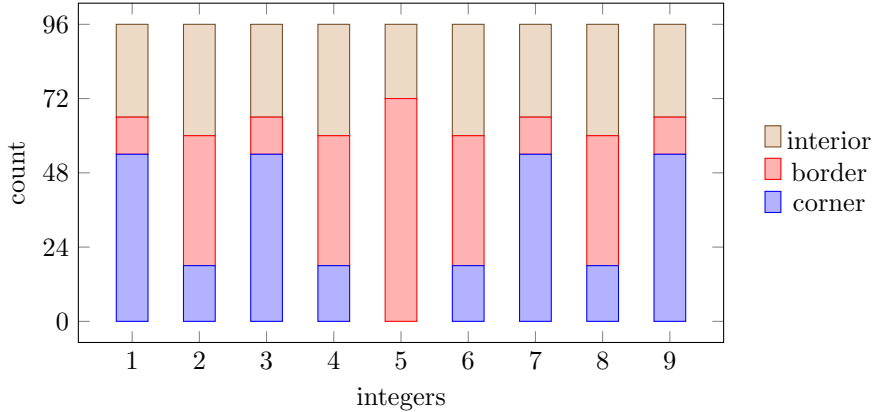
\begin{figure}[ht]
    \centering
    
    \begin{tikzpicture}
        \begin{axis}[width=4in, height=2.4in, ybar stacked, bar width=12pt, xtick=data, xlabel={integers},
        ylabel=count, ymin = 0, ymax = 96, ytick={0,24,48,72,96}, enlarge y limits = {abs=8pt}, reverse legend, legend style={at={(1.05,0.5)}, anchor=west, draw=none}]
        \addplot coordinates {(1,54) (2,18) (3,54) (4,18) (5,0) (6,18) (7,54) (8,18) (9,54)};
        \addplot coordinates {(1,12) (2,42) (3,12) (4,42) (5,72) (6,42) (7,12) (8,42) (9,12)};
        \addplot coordinates {(1,30) (2,36) (3,30) (4,36) (5,24) (6,36) (7,30) (8,36) (9,30)};
        \legend{corner, border, interior}
        \end{axis}
    \end{tikzpicture}
    
    \caption{Distribution of integers within all $3$-level magic triangles.}
    \label{fig:dist3}
\end{figure}

\Cref{fig:dist3} shows the distribution of each integer in all $3$-level magic triangles.
Intriguingly, the integer $5$ never appears in a corner entry, while it appears in a border entry in $72$ out of $96$ magic triangles. 
The other odd integers all have the same distribution, most often appearing in the corner entries and rarely appearing in border entries. 
Similarly, the even integers all have the same distribution, appearing least often in the corners but most often in the border entries.

\subsection*{4-level triangles}

The $4$-level triangle also contains orbits of entries under rotations and reflections. 
With indexing as in \cref{fourlevel}, we call $\{a_1, a_7, a_{16}\}$ the \emph{corner entries}, $A\cup B \cup C$ the \emph{border entries}, $\{a_2, a_6, a_{14}\}$ the \emph{interior entries}, and $a_{10}$ the \emph{center entry}.
We found that each integer $1, \ldots, 16$ occurs in each position with nearly equal frequency, though there are slight differences in the frequencies of integers in corner and border entries compared with interior and center entries.

\begin{figure}[ht]
    \centering
    
    \begin{tikzpicture}
        \begin{axis}[width=3.8in, height=2.4in, ybar stacked, bar width=8pt, xtick=data, xlabel={integers},
        ymin = 0, ymax = 1500, ytick={0,546,900,1174,1500}, yticklabels={0,$178{,}443{,}648$, $178{,}902{,}432$, $179{,}257{,}536$, $238{,}536{,}576$},
        ylabel=count, enlarge y limits = {abs=8pt}, reverse legend, legend style={at={(1.05,0.5)}, anchor=west, draw=none}]
            \addplot coordinates {(1,1002) (2,602) (3,1118) (4,546) (5,950) (6,986) (7,822) (8,1174) (9,1174) (10,822) (11,986) (12,950) (13,546) (14,1118) (15,602) (16,1002)};
            \addplot coordinates {(1,498) (2,898) (3,382) (4,954) (5,550) (6,514) (7,678) (8,326) (9,326) (10,678) (11,514) (12,550) (13,954) (14,382) (15,898) (16,498)};
            \legend{corner/border, interior/center};

        \end{axis}
        \begin{axis}[width=3.8in, height=2.4in, hide axis, xmin = 0, xmax = 16, ymin = 0, ymax = 1500, enlarge y limits = {abs=8pt}, clip=false]
            \addplot[dashed] coordinates {(0,900) (16,900)};
            \addplot[white, line width=4pt] coordinates {(-0.1,300) (16.1,300)};
            \addplot[white, line width=4pt] coordinates {(-0.1,1350) (16.1,1350)};
            \draw (-0.3,324)-- (0.3,324); 
            \draw (-0.3,276)-- (0.3,276);
            \draw (15.7,324)-- (16.3,324); 
            \draw (15.7,276)-- (16.3,276);
            \draw (-0.3,1374)-- (0.3,1374); 
            \draw (-0.3,1326)-- (0.3,1326);
            \draw (15.7,1374)-- (16.3,1374); 
            \draw (15.7,1326)-- (16.3,1326); 
        \end{axis}
    \end{tikzpicture}
    
    \caption{Distribution of integers within all $4$-level magic triangles.}
    \label{fig:dist4}
\end{figure}
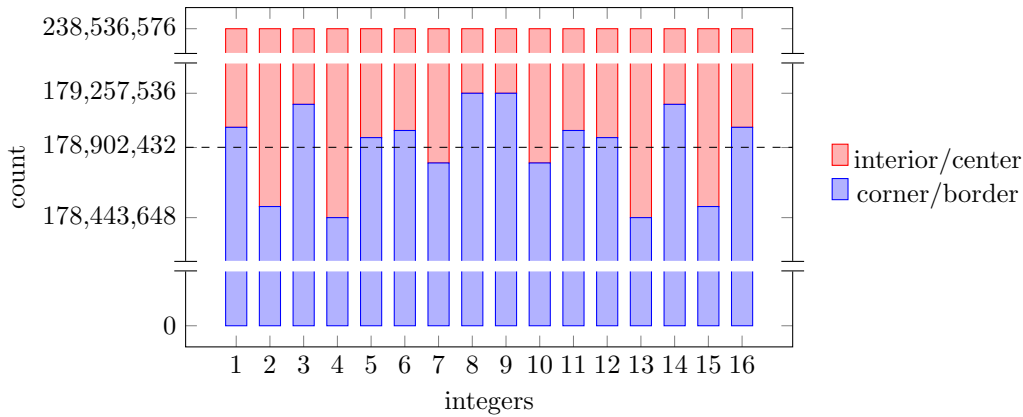

Surprisingly, any fixed $k \in \{1, \ldots, 16\}$ appears the same number of times in each corner and border position, and also the same number of times in each interior and center position.
However, if $k$ and $j$ are distinct positive integers such that $k+j \ne 17$, then $k$ and $j$ appear in each position with slightly different frequencies.
\cref{fig:dist4} displays the number of magic triangles in which each integer appears in the two combined groups corner/border and interior/center.
Note that there are $12$ corner/border entries and only $4$ interior/center entreis.
Thus, if an integer $k$ were equally distributed among all positions in all $4$-level magic triangles, then $k$ would appear in corner/border entries in exactly $\frac{3}{4}T_4 = 178{,}902{,}432$ magic triangles; the vertical axis in \cref{fig:dist4} has been constructed to highlight this value. 
While no integer occurs in exactly this many corner/border entries, \cref{fig:dist4} shows which integers occur slightly more often in corner/border entries, and which integers occur slightly more often in interior/center entries.

\section{Future Research} 

Our success at finding magic triangles using simulated annealing suggests that we could find an $n$-level magic triangle for any positive integer $n$, given enough computation time. 
However, we have not proved that magic triangles exist for all $n\in\mathbb{N}$.
We would like to find an algorithm for generating an $n$-level magic triangle for any $n \in \mathbb{N}$, which would prove that magic triangles of all sizes exist.

We have found the numbers of magic triangles for small $n$; the sequence $(T_n)$ of such numbers begins $1, 4, 96, 238536576$.
Finding the exact number of $n$-level magic triangles seems to be quite difficult for $n > 4$. 
We suspect that the number of magic triangles continues to grow quickly with $n$, though such triangles become less frequent among all triangular arrangements.
Perhaps randomized methods could be used to approximate the number of $n$-level magic triangles, as Kitajima and Kikuchi have done for magic squares \cite{Numerous}. 
Alternatively, it would be interesting to establish lower and/or upper bounds for the number of magic triangles for each $n$.

\section*{Code}

Python code for the computational experiments in this paper is available at\\
 \url{https://github.com/mlwright84/magictriangles}.

\section*{Acknowledgements}

This article grew out of the final project by the first two authors in \emph{Modern Computational Mathematics} at St.\ Olaf College in spring 2020.

\end{document}